\documentclass[11pt]{amsart}

\usepackage{epsfig}

\usepackage{amsfonts}
\usepackage{amssymb}
\usepackage{amsmath}
\usepackage{amsthm}
\usepackage{amscd}
\usepackage{amstext}
\usepackage{latexsym}

\usepackage{setspace}

\theoremstyle{plain}

\newtheorem{theorem}{Theorem}[section]

\newtheorem{proposition}[theorem]{Proposition}

\newtheorem*{division problem}{Division Problem}

\theoremstyle{definition}

\theoremstyle{remark}

\newcommand{\CC}{\mathbf{C}}

\newcommand{\RR}{\mathbf{R}}
\newcommand{\JJ}{\mathcal{J}}
\newcommand{\OO}{\mathcal{O}}

\newcommand{\noi}{\noindent}

\newcommand{\vp}{\varphi}
\newcommand{\qa}{\quad}

\newcommand{\lsm}{\preceq}

\providecommand{\abs}[1]{\left\vert #1\right \vert}
\newcommand{\xx}{\abs{x}^2}
\newcommand{\yy}{\abs{y}^2}
\newcommand{\zz}{\abs{z}^2}

\author{Dano Kim}

\title{A remark on the approximation of plurisubharmonic functions}

\begin{document}

\maketitle

\begin{abstract}

\noi We show by an example that the Demailly approximation sequence of a plurisubharmonic function, constructed via Bergman kernels, is not a decreasing sequence in general. 
  
\end{abstract}

\section{Introduction}

\label{}

 It is a fundamental theorem of Demailly~\cite{D92} that given an arbitrary plurisubharmonic function $\vp$ on a domain, there always exists a sequence $\{ \vp_m \}$ of plurisubharmonic functions with analytic singularities converging to $\vp$. Moreover, the approximating function $\vp_m$ is given in a very natural form: $\vp_m = \frac{1}{2m} \log \sum \abs{\sigma_l}^2 $ where $( \sigma_l )$ is an orthonormal basis of the Hilbert space of holomorphic functions that are square integrable with respect to the weight $e^{-2m \vp}$.  (See \cite[Proposition 3.1]{D92}, \cite{D09} and also the exposition in \cite{B}.)  

 It was further proved that the subsequence $\{ \vp_{2^k} \}$ is decreasing (after adding suitable constants) in \cite[Step 3, Proof of Theorem 2.3]{DPS} using  a subadditivity property of the sequence $\vp_m$'s. It remained a natural question, raised explicitly in \cite[p.134]{B}, to ask whether the entire sequence  $\{ \vp_m \}$ is decreasing or not, even up to the \emph{equivalence of singularities} (see the beginning of the next section for its definition). In this note, we show by an example that the Demailly approximation sequence of a plurisubharmonic function is not necessarily decreasing, thus answering the above question negatively.   The example $\vp$ is given as a plurisubharmonic function with analytic singularities, for which we can compute the multiplier ideal sheaf of each $m \vp$ and determine the singularities of $\vp_m$ using a finite number of local generators of $\JJ (m \vp)$.

\section{The Example} 

We recall the following definitions from \cite{D09}. For two singular weights  $h_i = e^{-\varphi_i} \; (i = 1,2)$, we say  $h_1$ is \textbf{less singular} than $h_2$ and write $h_1 \lsm h_2$ if the quotient $\frac{h_1}{h_2}$ is locally bounded above.  In this case, we also write $\vp_2 \lsm \vp_1$.  If $h_1 \lsm h_2$ and $h_2 \lsm h_1$, we say $h_1$ and $h_2$  \textbf{have equivalent singularities} and write $h_1 \sim h_2$ and $\vp_1 \sim \vp_2$.

 Let $X$ be a complex manifold and $D \subset X$ be an irreducible hypersurface. For a local defining function $f$ of $D$, we associate the plurisubharmonic function $\vp_D := \log \abs{f}^2$ (at least locally). More generally, let us denote by $\vp_D$ the plurisubharmonic function corresponding to an effective divisor $ D = \sum^{k}_{i=1} a_i D_i $, that is,  $\vp_D = \sum^k_{i=1} \log \abs{f_i}^{{2}{a_i}}$. Here $f_i$ is a local equation of $D_i$ and as is well-known, $\vp_D$ makes sense globally as local weight functions of a singular hermitian metric~\cite{D09} (of $\OO(D)$ if it is a line bundle). 
 
  We first give a trivial example of such $\vp_D$ when $D \subset X$ is a smooth hypersurface. In this case, $\JJ(m \vp_D) = \OO(-mD)$ which is locally generated by one generator $f^m$ where $D = (f=0)$. It follows that in this case, the Demailly approximation $\{ \vp_m \}$ is a constant sequence, at least up to equivalence of singularities. 

 For a  simple example with nonconstant $\{ \vp_m \}$, let $X = \CC^2$  (with coordinates $x,y$)  and $\vp = \log (\abs{x}^2 + \abs{y}^2)$. Since $\JJ(m \vp) = \mathfrak{m}^{m-1}$ (where $\mathfrak{m} \subset \OO_X$ is the maximal ideal of the point $(0,0)$),  it follows from considering the generators of the ideal $\mathfrak{m}^{m-1}$  that $\displaystyle \vp_m = \frac{1}{m} \log ( \sum^{m-1}_{i=0} \abs{x^{i} y^{m-1-i} }^2 ) \sim \frac{m-1}{m} \log (\abs{x}^2  + \abs{y}^2) $, which is strictly decreasing in terms of equivalence of singularities.

 In order to construct our example, let us take $X = \CC^2$ (with coordinates $x,y$) and the effective divisor $D = \sum^{3}_{i=1} a_i D_i  \qa (a_i \ge 0)$, where $D_1 = (x=0)$, $D_2 = (y=0)$, $D_3 = (z:= x+y =0)$.  This is the simplest example of a non-SNC divisor. One blowup of the origin $\pi: X' \to X$ gives the log-resolution of the pair $(X,D)$ and we can compute the multiplier ideal sheaves.
 
 Let $H_i \subset X'$ be the proper transform of $D_i$ and $E$ the exceptional divisor of $\pi$. Then the multiplier ideal sheaves are given by  (\cite{D09}, \cite{L})

$$  \JJ(m \vp_D) = \JJ(mD) = \pi_* \OO_{X'} ( K_{X'} - \pi^* (K_X) - \lfloor \pi^* (mD) \rfloor  )   $$
   $$  = \pi_* \OO_{X'} ( - \sum^3_{i=1} \lfloor m a_i \rfloor  H_i   - (2m -1) E). $$

 Now let $\Omega \subset \CC^2$ be a connected Stein neighborhood of the origin and  consider the Demailly approximation sequence $\vp_m$ on $\Omega$  \cite[Proposition 3.1]{D92} given by the Bergman kernels  corresponding to the weight function $e^{- m \vp_D}$. From now on, we identify a plurisubharmonic function with its equivalence class in terms of singularities as in \cite[Definition 6.3]{D09}. If the sequence $\{ \vp_m \}$ is genuinely decreasing (up to some constants), then of course it is also decreasing in terms of equivalence of singularities. 
 
 Since $\vp_m$ is equivalent to a plurisubharmonic function with analytic singularities given by a finite number of generators of the multiplier ideal sheaf $\JJ (m \vp_D)$ in a relatively compact Stein open subset of $\Omega$, we can use it to show the following. 
\\

\begin{theorem}\label{exx}

 Let $\{ \vp_m \}$ be the Demailly approximation sequence for the plurisubharmonic function $\vp = \vp_D$ where $D$ is as in the above with $a_1 = a_2 = a_3 = \frac{2}{3}$. Then the sequence is not decreasing in the sense that we cannot choose a sequence of constants $C_m$ such that $\vp_m + C_m$ is a decreasing sequence of $\RR \cup \{ - \infty \}$-valued functions. 

\end{theorem}
\quad
\\

The coefficients $a_1 = a_2 = a_3 = \frac{2}{3}$ are chosen here only for the simplicity of computation. Apparently, the behavior of the approximation sequence as in this theorem are expected to be very common for other psh functions. 


\noindent {\it{Proof of Theorem \ref{exx}}.}   Let $\vp = \vp_D$. From now on,  functions have the domain as a unit ball around the origin. For $m \ge 1$, the psh function $2 \vp_m$ is equivalent to $\frac{1}{m} \log (\sum \abs{f_i}^2 )$ where $f_i$'s are the finite number of generators for $\JJ(m \vp)$.  As for the multiplier ideal sheaf $\JJ(2\vp) = \pi_* \OO_{X'} (- \sum H'_i - 3E)$, we see that  
it can be generated (over $\OO_X$) by $f_1 = xyz$. Hence for $m=2$, we have $ 2 \vp_2 \sim \log \abs{xyz}$. Similarly, for $\JJ(3\vp) = \pi_* \OO_{X'} (- \sum 2 H'_i - 5E)$, the generator can be chosen as $(xyz)^2$, thus $ 2 \vp_3 \sim \log \abs{xyz}^{\frac{4}{3}}$.

 One then immediately sees that the sequence of singularities is not decreasing:   $ \vp_5 \npreceq \vp_3 $  since $2 \vp_5 \sim \log \abs{xyz}^{\frac{6}{5}}$. In fact, we see that $\vp_5 \succeq \vp_3$ holds in this case. Also one checks that $\vp_4 \npreceq \vp_3$:  ($\JJ (4 \vp)$ has generators $(xyz)^2 x$, $(xyz)^2 y$, $(xyz)^2 z$, the last being redundant but included for convenience)

 $$ \frac{e^{2 \vp_4}}{e^{2 \vp_3}} \sim  \frac{  \abs{xyz} (\xx + \yy + \zz)^{\frac{1}{4}}   }{ \abs{xyz}^{\frac{4}{3}}    }\qa   $$  is  not locally bounded above, considering along $x=y$.  This completes the proof of Theorem~\ref{exx}. 
 
 Furthermore, similar computations yield that there is an infinite number of instances where the decreasing property fails to hold, considering $m$ modulo $3$ for the coefficient $\lfloor \frac{2m}{3} \rfloor $. We have
 
$$ 2 \vp_{3k} \sim \log \abs{xyz}^{2 \frac{2k}{3k}}   \preceq    2 \vp_{3k+2} \sim \log \abs{xyz}^{2 \frac{2k+1}{3k+2}}  $$ generalizing $\vp_5 \succeq \vp_3$. So we cannot truncate the sequence $\{ \vp_m \}$ to make it decreasing.


\quad
\\

 On the other hand, similar considerations can be used to check the decreasing property of the  subsequence with indices of exponential growth \cite[Step 3, Proof of Theorem 2.3]{DPS} of our main example in Theorem~\ref{exx} directly, without using subadditivity:
\\

 \begin{proposition}\label{dec}
  $ \vp_{2^{k+1}} \preceq \vp_{2^k} $ for every $k \ge 1$.  
  
 \end{proposition}

\noindent For the proof of this,  we use the facts that $\lfloor 2^{2k} \frac{2}{3} \rfloor = 2 \frac{2^{2k} -1}{3}$, $\lfloor 2^{2k+1} \frac{2}{3} \rfloor = 2 \frac{2^{2k+2} -1}{3}$ and that $ \xx + \yy + \zz   \ge  3 \abs{xyz}^{\frac{2}{3}} $. For the first few terms, we find
 
 $\displaystyle \frac{e^{2 \vp_8}}{e^{2 \vp_4}} \sim  \frac{ \abs{xyz}^{\frac{10}{8}}    }{  \abs{xyz} (\xx + \yy + \zz)^{\frac{1}{4}} } \preceq \abs{xyz}^{\frac{1}{2}} $    is locally bounded above.

 $\displaystyle \frac{e^{2 \vp_{16}}}{e^{2 \vp_8}}  \sim  \frac{ \abs{xyz}^{\frac{10}{8}} (\xx + \yy + \zz)^{\frac{1}{16} } }  {\abs{xyz}^{\frac{10}{8}} }  $ is locally bounded above thanks to the equality of the exponents of $\abs{xyz}$ in the fraction. Then it is easy to check that  these two patterns alternate and show $ \vp_{2^{k+1}} \preceq \vp_{2^k} $ for every $k \ge 1$. This completes the proof of Proposition~\ref{dec}. 
\\
  
  Finally, as Professor J.-P. Demailly kindly pointed out to us, we note that it might still be possible to show that there exists a (strictly) decreasing subsequence of the Demailly approximation with indices of linear growth, instead of exponential growth as above. In our main example, this is indeed the case since the subsequence $\vp_{3k+2} \sim  \log \abs{xyz}^{ \frac{2k+1}{3k+2}} $ is strictly decreasing and does converge to the original $\vp$. The arguments used in this note might be possibly used for general $\vp$ with analytic singularities. 





\section*{Acknowledgements}

This work was supported by the  National Research Foundation of Korea grants NRF-2012R1A1A1042764 and No.2011-0030795, funded by the Korea government and also by Research Startup Fund for new faculty of Seoul National University.

\qa

\qa

\normalsize

\noi \textsc{Dano Kim}

\noi Department of Mathematical Sciences, Seoul National University

\noi 1 Kwanak-ro, Kwanak-gu, Seoul, Korea 151-747

\noi Email address: kimdano@snu.ac.kr

\noi

\end{document}